\documentclass[12pt]{article}
\usepackage{graphicx} 
\usepackage{amsmath,amssymb}
\usepackage{amsthm}   
\theoremstyle{remark}

\usepackage{amsfonts}
\usepackage{amsmath}
\usepackage{multirow} 
\usepackage{booktabs}
\usepackage{pdflscape}
\usepackage{float}
\usepackage{amssymb}
\usepackage{algorithm}
\usepackage{algorithmic}
\usepackage[numbers]{natbib}
\usepackage{hyperref}
\usepackage{textcomp}
\usepackage{xcolor}
\usepackage{makecell}
\definecolor{highlight}{RGB}{133,223,255} 
\usepackage[dvipsnames,table]{xcolor}
\usepackage{booktabs} 
\usepackage[a4paper, total={6.8in, 9in}]{geometry}
\title{A hybrid wavelet-based physics-informed neural network  for portfolio management}

\author{
Bahadur Yadav$^{1}$,
Mahaprasad Mohanty$^{1}$,\\
Ratikanta Behera$^{2}$,
Sanjay Kumar Mohanty$^{1,*}$
}

\date{}

\begin{document}

\maketitle

\begin{center}
$^{1}$Department of Mathematics, School of Advanced Sciences,\\
Vellore Institute of Technology, Vellore, Tamil Nadu 632014, India\\[6pt]

$^{2}$Department of Computational and Data Sciences,\\
Indian Institute of Science, Bangalore 560012, India \\[6pt]
$^{*}$Corresponding author: \href{mailto:sanjaykumar.mohanty@vit.ac.in}{sanjaykumar.mohanty@vit.ac.in}
\end{center}

\abstract{In this paper, we present a Hybrid Wavelet-based Physics-Informed Neural Networks (HW-PINNs) framework for portfolio management that provides a promising alternative to Physics-Informed Neural Networks (PINNs). Here, we first discuss the generalized framework of the Merton jump diffusion model and the associated HW-PINNs, followed by the one-dimensional case of a European option. Our work adapts the HW-PINN framework to the Merton jump-diffusion model for a European option, using a simplified direct coefficient optimization strategy, a mathematically corrected log-space formulation, and an efficient FFT -based computation of the integro-differential operator. Through numerical experiments across realistic market scenarios, we show that our proposed model achieves high accuracy and robustness, with a mean relative error of 0.27\% in low jump intensity scenarios compared to high-fidelity benchmarks. Our results validate that the implementation of this specific HW-PINN framework is a computationally efficient and reliable tool for pricing derivatives in markets with high jump risk. In addition, we further discuss risk analysis using Value at Risk (VaR) and Conditional Value at Risk (CVaR), which provide insights into downside risk across different market scenarios.}
 
\noindent\textbf{Keywords:} {Physics-Informed Neural Networks; Wavelets; Merton Jump-Diffusion Model; European Option Pricing; Fast Fourier Transform; Risk Analysis.}


\maketitle

\section{Introduction}
Portfolio optimization is a branch of financial mathematics that has gained significant popularity over the past several decades due to its practical relevance and wide-ranging real-life applications. It was initially established in frictionless, return-based environments that cannot account for tail risks, nonlinear payoffs, or optimistic market expectations. It is one of the traditional mathematical finance challenges that has been thoroughly studied for the simultaneous maximization of returns and minimizing of risk. Early studies have shown that realistic elements, such as constant and proportional transaction costs, strongly influence optimal investment behavior by creating no-trade regions. For derivative-based methods, this phenomenon is particularly significant \cite{oksendal2002optimal}. It was further reframed as a signal-extraction and convex-optimization problem, providing reliable methods for handling high-dimensional data and estimation error \cite{zhao2019optimal}.
A significant transition occurred when it was demonstrated that option prices incorporated forward-looking information, with option-implied volatility and asymmetry boosting portfolio performance relative to historical estimates \cite{demiguel2013improving}. This realization led to the direct optimization of option-based portfolios, specifically addressing asymmetric and non-Gaussian reward structures that extend beyond mean-variance analysis \cite{faias2017optimal}.
Later, dynamic trading techniques were introduced that enhance intertemporal portfolio efficiency by incorporating derivatives to hedge against jump risk and stochastic volatility \cite{liu2003dynamic}. Although they were more sensitive to option surface estimate, contemporary approaches extracted complete option-implied distributions and state variables, allowing portfolio decisions that reflect market-implied tail risk \cite{kostakis2011market}. However, further researchers developed a technique to manage skewness, kurtosis, and nonlinear exposures using option-based instruments due to the significance of higher-order moments in portfolio risk \cite{bhandari2009options}. To highlight hypergeometric stochastic volatility as a distinct and tradable risk factor, it can be strategically managed and hedged through the use of option-based instruments \cite{cipriano2021optimal}.

 Also, option prices offer forward-looking information about volatility and tail risk. Jump-diffusion frameworks are used because continuous diffusion models cannot capture rapid market fluctuations. To capture such phenomena, the jump–diffusion framework was introduced, in which asset prices evolve through a combination of continuous Brownian motion and discrete jumps \cite{merton1976option}. This method gives a more accurate depiction of asset dynamics and serves as the foundation for expanding portfolio optimization challenges where discontinuities are important. Further studies have examined how the presence of jumps affects optimal investing strategies, showing that ignoring jumps may result in less than ideal allocation choices and that optimal portfolio weights can be broken down into diffusion-driven and jump-risk features \cite{jin2012decomposition}. Later, research studies have examined the circumstances in which jumps have a substantial impact on ideal portfolios, demonstrating that the significance of jump risk is dependent on variables such as investor risk choices, jump frequency, and jump size distribution \cite{ascheberg2016jumps}. Furthermore, risk-sensitive portfolio optimization methods have been created to specifically account for regime shifts and jump risks, enabling investors to modify their strategies in response to shifting market conditions and degrees of uncertainty \cite{das2018risk, michelbrink2012martingale}. In order to optimize portfolios in more realistic multi-asset markets, the jump-diffusion technique was extended to multivariate situations where several assets undergo corresponding jump circumstances \cite{de2025portfolio}. 

European option–based portfolio optimization can be viewed as a special case of a higher-dimensional framework, where investment decisions are tied to derivative payoffs written on a single underlying asset. Initial studies on European options under market imperfections examined pricing and hedging in the presence of fixed and proportional transaction costs, showing that trading frictions significantly affect optimal hedging strategies and portfolio choices \cite{zakamouline2006european}. Furthermore, research studies on derivative-based investment strategies from an expected utility perspective have demonstrated that European options in portfolios can improve overall welfare and risk management by exposing investors to nonlinear payoff structures \cite{escobar2022derivatives}. Financial markets are often incomplete and exposed to credit risk, and subsequent research has shown how default risk and market imperfections can significantly influence both portfolio pricing and optimal allocation decisions \cite{grigorova2020european}. In order to capture market-implied dynamics, nonparametric calibration procedures were developed to analyze jump-diffusion option pricing models from a practical modeling approach \cite{cont2004nonparametric}. Advanced portfolio selection strategies under uncertainty that account for competing investment objectives and model uncertainty in portfolio formation have been the focus of recent studies \cite{mohseny2025robust, yadav2025hybrid}. For adaptive portfolio models that incorporate risk-free return perspectives, option strategies, and shadow pricing methods have been developed to enhance decision-making in complicated financial environments \cite{escobar2022derivatives}. 

Conventional methods depend on numerical techniques, such as Fourier methods and Monte Carlo simulation, to solve option pricing problems under jump-diffusion dynamics. When applied to high-dimensional and nonlinear financial frameworks, these approaches often face computational complexity, grid dependence, and dimensionality costs, despite their theoretical robustness. Physics-Informed Neural Networks (PINNs) were initially introduced within machine learning as a way to incorporate governing financial PDEs directly into the training process, enabling efficient approximation of solutions to differential equations \cite{raissi2019physics}. A variety of differential equations that arise have been solved using this framework \cite{shi2024physics}.
 Furthermore, deep learning-based algorithms have shown favorable outcomes in resolving the curse of dimensionality when solving high-dimensional PDEs \cite{krishnapriyan2021characterizing}. These advancements demonstrate that neural network approaches can effectively approximate solutions to sophisticated PDE systems without relying on conventional mesh-based discretization methods \cite{sirignano2018dgm}. However, more sophisticated neural network frameworks, such as deep learning techniques for volatility, hedging, and derivative pricing, are being discussed for financial markets \cite{buehler2019deep,  beck2021deep}. Moreover, neural network approaches to parabolic PDEs have improved numerical reliability and efficiency for high-dimensional problems \cite{horvath2021deep}. Additionally, recent studies have improved derivative pricing accuracy by combining deep learning with asymptotic methods \cite{funahashi2026deep}. 
    However, jump-diffusion models often have sharp gradients and discontinuities because of abrupt price jumps and integral components in the governing Partial Integro-Differential Equation (PIDE). In these situations, ordinary PINNs are not efficient. 
 
 To deal with this specific challenge, the Wavelet-based Physics-Informed Neural Networks architecture has been introduced as a powerful tool \cite{pandey2024efficient}. In our framework, we discuss the Merton jump-diffusion model for European options to optimize portfolios. Our proposed approach combines the computational efficiency of the Fast Fourier Transform with a wavelet-based spectral approach to primarily address these shortcomings. Using multi-resolution analysis, HW-PINNs represent a solution as a linear combination of wavelet family functions rather than as the direct output of a neural network.  The most effective coefficients for this wavelet series are found by training the network. Our approach is especially useful for two reasons: first, wavelets are well-suited to describing localized, multi-scale phenomena that are difficult for regular PINNs to capture due to their strong time-frequency localization. Second, the derivatives needed for the PDE residual can be calculated analytically by building the solution from an analytical foundation, avoiding automatic differentiation and the processing overhead that goes along with it. Comparing our FFT-based approach to conventional quadrature demonstrates two significant advantages. First, the FFT algorithm's computational complexity is $O(N \log N)$, which is significantly faster than the $O(N^2)$ complexity of direct numerical integration for $N$ points. This is the most noticeable aspect. Second, our technique prevents the compounding numerical inaccuracies that can result from repeated, localized interpolations within a naive integration loop by working on the entire solution vector at once. We combine FFT-based computation of the jump integral with the analytical derivatives from the wavelet matrices. Furthermore, our results are often very similar to the benchmark values, demonstrating the model's better accuracy under various market circumstances. In particular, our model provides better results with fewer errors in the low jump intensity case. In order to assess downside risk under various market conditions, we further discuss a risk analysis using VaR and CVaR.

This paper is organized as follows: We begin by giving the mathematical background of our proposed approach in Section \ref{section:2}.
Section \ref{section:3} provides the problem formulation of our proposed approach. In Section \ref{section:4}, we discuss HW-PINN architectures, techniques for incorporating physical laws, and the solution of the Merton PIDE. Section \ref{section:5} provides empirical results that support the effectiveness of our approach. Finally, Section \ref{section:6} concludes the paper and outlines potential directions for future research.

\section{Portfolio Optimization under Jump Risk}\label{section:2}
 In this section, we discuss a general Merton jump-diffusion model, which provides a unified framework for asset prices that ensures positivity of asset prices and analyze the portfolio optimization by leveraging the option price. 
Considered the Merton jump-diffusion model on a filtered probability space
\(\
(\Omega,\mathcal{F},\mathbb{P},(\mathcal{F}_t)_{t\ge 0}),
\)
where $\Omega$ represents the sample space, $\mathcal{F}$ denotes a $\sigma$-algebra on $\Omega$, $\mathbb{P}$ represents the probability measure, and $(\mathcal{F}_t)_{t\ge 0}$ is a right-continuous filtration that models the flow of information available to market participants over time \cite{merton1976option}. This model determines the admissible dependence of asset prices and trading strategies on past and present information. In the Merton jump--diffusion model, asset prices are influenced by two different kinds of uncertainty: continuous market fluctuations and sudden discontinuous jumps. To capture these features, the underlying sample space is constructed as the product $\Omega=\Omega_W\times\Omega_N$, where $\Omega_W=C([0,\infty),\mathbb{R}^d)$ represents the space of continuous Brownian paths, modeling the diffusion component of asset returns, and $\Omega_N$ represents the space of jump counting paths associated with a Poisson process.  Now, the probability measure $\mathbb{P}$ on $\Omega$ is defined as the product measure $\mathbb{P}=\mathbb{P}_W\otimes\mathbb{P}_N$, which ensures the independence between the continuous diffusion component and the jump component of the market dynamics. In the probability space $(\Omega_W,\mathcal{F}_W,\mathbb{P}_W)$, the process $W=(W_t)_{t\ge0}$ is a $d$--dimensional Brownian motion, where $d\in\mathbb{N}$. The process $W_t$ has independent and stationary increments that satisfy $W_t-W_s\sim\mathcal{N}(0,(t-s)I_d)$. However, in the probability space $(\Omega_N,\mathcal{F}_N,\mathbb{P}_N)$, the process $N_t$ is a Poisson process with intensity $\lambda>0$, satisfying
\[\
\mathbb{P}(N_t=k)=\frac{(\lambda t)^k e^{-\lambda t}}{k!}.
\]

Furthermore, jump sizes $(Y_i)_{i\ge1}$ are independent and identically distributed and assumed to be independent of both the Brownian motion $W$ and the Poisson process $N$. The associated compound Poisson process is written as
\(\
J_t=\sum_{i=1}^{N_t}Y_i
\)
represents the cumulative impact of jumps on the asset price up to time $t$. Then, the market filtration is generated by the observable sources of randomness, and is given by
\(\
\mathcal{F}_t=\sigma\big(W_s,N_s,Y_i:\, s\le t\big),
\)
and is enhanced in the standard manner. Therefore, any financially admissible asset price process must be adapted to $\mathcal{F}_t$,    guaranteeing that prices at time $t$ rely only on data that has been available up to that time.

Now, we consider \(\ R_t = \big(R_t^{(i)}\big)_{i=1,\dots,d}
\) be a $d$-dimensional Merton jump--diffusion process that models the logarithmic returns of the risky assets \cite{merton1976option}. The return process is given as
\[
R_t
=
\mu t  
+
\Sigma^{1/2} W_t
+
\sum_{k=1}^{N_t} Y_k,
\qquad t \ge 0,
\]
where $\mu\in\mathbb{R}$ is the drift parameter of the return process, $\Sigma^{1/2}$ is the square root of the diffusion covariance matrix, $W_t$ is the Brownian motion defined on probbility space $(\Omega_W,\mathcal{F}_W,\mathbb{P}_W)$, and $\sum_{k=1}^{N_t} Y_k$ denotes the cumulative jump contribution arising from the compound Poisson process.

Further, we consider discounted asset prices defined as the component-wise exponential transformation, and written as
\(\
S_t
=
\big(
S_t^{(1)},\dots,S_t^{(d)}
\big)
=
\big(
S_0^{(1)} e^{R_t^{(1)}},\dots,S_0^{(d)} e^{R_t^{(d)}}
\big),
\ t \ge 0.
\label{eq:asset_prices}
\)
We also include a risk-free asset given by
\(\
S_t^{(0)} \equiv 1
\) . This specification ensures that asset prices remain strictly positive.
The price process follows the stochastic exponential representation and is defined as
\(\
\frac{S_t}{S_0}
=
e^{R_t}
=
\mathcal C(\hat R)_t,
 \ t \ge 0,
\)
where $\mathcal C(\cdot)$ represents the stochastic exponential and $\hat R_t$ denotes the semi-martingale driving the gain process, whose jumps satisfy
\(\
\Delta \hat R_t = e^{\Delta R_t}-1.
\)
We denote jumps of the return process by
\(\
\Delta R_t = J,\
\Delta \hat R_t = \hat J.
\)

Let \(\
V:\mathbb R^d \times [0,T] \to \mathbb R \) represents the value function of the European option contingent whose payoff depends on a $d$-dimensional vector of volatile assets. It is assumed that \(V\) is sufficiently smooth, that is,  \(V \in C^{2,1}\), for the multidimensional Itô formula to be applicable \cite{oksendal2007applied}. In the process \(V(S_t,t)\), we get
\[\
dV
=
V_t\,dt
+ \nabla_S V^\top\, dS_t
+ \tfrac12\, dS_t^\top H_V\, dS_t 
+ \big[V(YS_t,t)-V(S_t,t)\big]\, dN_t.
\]
Here, \(\nabla_S V\) represents the gradient of \(V\) with respect to \(S\)
and \(H_V\) is its Hessian matrix. The discontinuous change in the option value driven by a jump in the underlying asset values is captured by the last term. 
This expression describes the infinitesimal evolution of the option value along the stochastic dynamics of the asset price process \(S_t\).

Now, the risk-free probability measure requires that the discounted option price be a martingale. Therefore, taking expectations eliminates the stochastic components and results in a deterministic evolution equation of the following form:
\[\
\partial_t V + \mathcal{A}V = 0 ,
\] where \(\mathcal{A}\) represents the infinitesimal generator of the jump diffusion process. The operator \(\mathcal{A}\) contains the diffusion component from the Brownian motion and the nonlocal jump component coming from the Poisson process.
For numerical computation, it is convenient to define log-price variables as follows:
\(\
x=\log S, 
\
u(x,t)=V(S,t),
\)
which transforms the pricing equation into a linear partial integro-differential equation and is given as
\[\
\partial_t u = \mathcal{V}u ,
\]
where \(\mathcal{V}\) represents a L\'evy operator mapping on functions
\(u:\mathbb{R}^d \to \mathbb{R}\). Thus, the Fourier transform diagonalizes L\'evy  operators. Therefore, by applying the Fourier transform, we get
\[\
\hat{u}(\eta,t) = \mathcal{F}[u(x,t)],
\] and the equation becomes
\[\
\partial_t \hat{u}(\eta,t) = \Phi(\eta)\hat{u}(\eta,t).
\]
Here \(\Phi(\eta)\) is a characteristic exponent of the
\(d\)-dimensional jump-diffusion process.

Solving this ordinary differential equation, we yield
\[\
\hat{u}(\eta,t) = e^{t\Phi(\eta)}\hat{u}_0(\eta).
\] Now, taking the inverse Fourier transform, we obtain 
\[\
u(x,t)
=
\mathcal{F}^{-1}
\left(
e^{t\Phi(\eta)}
\mathcal{F}(u_0)
\right).
\] Discretizing the spatial domain and replacing the continuous Fourier transform with the discrete Fourier transform goes directly to the FFT approximation
\[\
u_N(x,t)
=
\mathcal{F}_N^{-1}
\left(
e^{t\Phi(\eta_k)}
\mathcal{F}_N(u_0)
\right).
\]
As a result, the jump-diffusion PIDE can be solved numerically efficiently using the FFT. 

We further examine a financial market in which a multivariate variation of the Merton jump diffusion model governs the underlying asset values. Let the vector of asset prices be represented by $S = (S_1,\dots,S_d) \in \mathbb{R}^d$. A partial integro-differential equation of the following form is satisfied by the associated option price $u(x,t)$ \cite{schoutens2003levy} and written as:
\[\
\frac{\partial u}{\partial t}
+
\mathcal{L}u
+
\mathcal{J} u
= 0,
\quad x \in \mathbb{R}^d.
\]
The diffusion and jump operators, which are derived using Itô's lemma and expectation under the risk-neutral measure, incorporate the stochastic dynamics of the underlying assets into the deterministic PIDE. In this case, the continuous diffusion dynamics are captured by the differential operator $\mathcal{L}$, which is given by:
\[\
\mathcal{L} u
=
\sum_{i=1}^{d}
\mu_i S_i \frac{\partial u}{\partial S_i}
+
\frac{1}{2}
\sum_{i=1}^{d}
\sigma_i^2 S_i^2 \frac{\partial^2 u}{\partial S_i^2},
\]
where the drift and volatility of the $i$-th asset are represented by $\mu_i$ and $\sigma_i$, respectively.
The nonlocal operator describes the jump component and is given by:
\[\
\mathcal{J} u
=
\int_{\mathbb{R}^d}
\big[
u(x+y,t) - u(x,t)
\big] \nu(dy),
\]
where $\nu$ represent the Lévy measure characterizing the jump distribution.

In general, Fourier basis functions are used to represent the solution in the frequency domain in FFT-based techniques. Nevertheless, spatial localization does not occur from these representations. In order to solve this, we use a wavelet-based representation that incorporates multiscale and localized basis functions to generalize the Fourier expansion \cite{ernst2025certified, uddin2023wavelets}. We employ a multiscale wavelet basis to approximate the solution. In particular, the solution is given as:
\[\
u_\theta(x,t)
=
\sum_{j=j_0}^{J}
\sum_{k \in \mathbb{Z}^d}
c_{j,k}(t)\, \Phi_{j,k}(x),
\quad x \in \mathbb{R}^d,
\]
where $c_{j,k}(t)$ are the random coefficients and $\theta$ represents the set of trainable parameters.
A family of $d$-dimensional wavelet basis functions is formed by the functions $\Phi_{j,k}(x)$, which are defined by:
\[\
\Phi_{j,k}(x)
=
2^{\frac{jd}{2}} \Phi(2^j x - k),
\]
which guarantees suitable normalization in $L^2(\mathbb{R}^d)$. 
A tensor-product construction is often used for practical implementation and is given by:
\[\
\Phi_{j,k}(x_1,\dots,x_d)
=
2^{\frac{jd}{2}}
\prod_{i=1}^{d}
\Phi(2^j x_i - k_i).
\]
Our model can capture both global structures and localized aspects of the solution because of its multiscale representation, which is crucial when there are jumps and non-smooth payoffs. 

In order to integrate the governing PIDE into the learning framework, we substitute the equation with the wavelet-based approximation $u_\theta(x,t)$. This gives a multiscale approximation of the solution \cite{ernst2025certified, han2025physics}.  This leads to the following physics-informed formulation:
\[\
\frac{\partial u_\theta}{\partial t}
+
\mathcal{L} u_\theta
+
\mathcal{J} u_\theta
= 0.
\]

We further define the corresponding residual function as follows:
\[\
\mathcal{R}_\theta(x,t)
=
\frac{\partial u_{\theta}}{\partial t}
+
\mathcal{L} u_{\theta}
+
\mathcal{J} u_\theta.
\]

In order to effectively enforce the controlling PIDE at certain collocation points, the HW-PINN aims to minimize this residual across the domain.
By minimizing a composite loss function that enforces the PIDE and the corresponding initial and boundary conditions, the parameters $\theta$ are determined as follow:
\[\
\mathcal{L}_{\text{HW-PINN}}
=
\frac{1}{N_c}
\sum_{i=1}^{N_c}
\left|
\mathcal{R}_\theta(x_i,t_i)
\right|^2
+
\mathcal{L}_{\text{IC}}
+
\mathcal{L}_{\text{BC}},
\]
The collocation points in the computational domain are denoted by $\{(x_i,t_i)\}_{i=1}^{N_c}$. The initial condition loss is given by:
\[\
\mathcal{L}_{\text{IC}}
=
\frac{1}{N_{ic}}
\sum_{i=1}^{N_{ic}}
\left|
u_\theta(x_i,0) - u_0(x_i)
\right|^2,
\]
where the payoff function in the financial sector is represented by $u_0(x)$, which is the specified initial condition.  The boundary condition loss is defined as:
\[\
\mathcal{L}_{\text{BC}}
=
\frac{1}{N_{bc}}
\sum_{i=1}^{N_{bc}}
\left|
u_\theta(x_i,t_i) - h(x_i,t_i)
\right|^2,
\]
where the asymptotic characteristics of the option price determine the boundary condition of the solution, which is represented by $h(x,t)$.

In the following section, we will discuss the one-dimensional case of the option pricing for a portfolio optimization in a jump–diffusion market.

\section{European Option in the Merton Jump Diffusion Model}\label{section:3}
In this section, we discuss the European option in the one-dimensional case of the Merton jump-diffusion model for portfolio optimization. Let
\(
S_t = (S_t^1,\dots,S_t^d)
\)
represent the vector of prices of $d$ volatile assets. Each asset price follows the jump--diffusion
dynamics, and is given as
\begin{equation*}
\frac{dS_t^i}{S_{t-}^i}
=
\mu_i\,dt
+
\sum_{j=1}^d \sigma_{ij}\,dW_t^j
+
(Y_i-1)\,dN_t,
\qquad i=1,\dots,d,
\label{eq:d_dim_merton}
\end{equation*}

\noindent where $\mu_i$ denotes the drift of the  $i$-th asset, $\sigma_{ij}$ represents its exposure to the $j$-th
Brownian source, and $Y_i>0$ is a jump multiplier. At a jump time $t$, the price of the asset
is given by $S_t^i = S_{t-}^i Y_i$, ensuring the positivity of the prices.

Despite being presented in a $d$-dimensional environment, the pricing of a European option expressed on a single asset is only determined by the asset's marginal dynamics. Assuming that $S_t^i$ is fixed, the diffusion term \(\sum_{j=1}^d \sigma_{ij}\,dW_t^j \) is Gaussian with instantaneous variance \(\sum_{j=1}^d \sigma_{ij}^2\,dt \). Further, when the effective volatility is defined as follows: \(
\sigma_i^2 := \sum_{j=1}^d \sigma_{ij}^2
\), the marginal dynamics reduce to the one-dimensional Merton jump--diffusion model is given as
\begin{equation}
\frac{dS_t}{S_{t-}}
=
\mu\,dt
+
\sigma\,dW_t
+
(Y-1)\,dN_t.
\label{eq:1d_merton}
\end{equation}

Now, let
\(
V:\mathbb R^+ \times [0,T] \to \mathbb R
\)
is a value function of a European option defined on the asset $S_t$, and
assume that $V \in C^{2,1}$. By applying the It\^o--L\'evy formula to $V(S_t,t)$, we get
\begin{align}
dV
&=
V_t\,dt
+
V_S\,dS_t
+
\frac12 \sigma^2 S_t^2 V_{SS}\,dt
+
\big[ V(YS_{t^-},t) - V(S_{t^-},t) \big]\,dN_t.
\label{eq:2d_merton}
\end{align}
Substituting the asset dynamics \eqref{eq:1d_merton} in \eqref{eq:2d_merton}, gives the following
\begin{align}
dV
&=
\Big(
V_t
+
\mu S V_S
+
\frac12 \sigma^2 S^2 V_{SS}
\Big)dt
+
\sigma S V_S\,dW_t
+
\big[ V(YS,t) - V(S,t) \big]\,dN_t.
\end{align}

The drift term $\mu$ is replaced by the risk-free interest rate $r$ in an equivalent martingale measure $\mathbb Q$. When the no-arbitrage condition is satisfied, then  the partial integro-differential equation (PIDE) for the value of a European option is as follows:

\begin{equation}
\frac{\partial V}{\partial t}
+ \frac12 \sigma^2 S^2 \frac{\partial^2 V}{\partial S^2}
+ r S \frac{\partial V}{\partial S}
- (r+\lambda) V
+ \lambda \int_{-\infty}^{\infty} V(Se^y,t)\,f_Y(y)\,dy
= 0.
\label{eq:merton_pide}
\end{equation}

Here, $f_Y$ denotes the probability density function of the logarithmic jump size. The integral term captures the contribution of jump risk, while the differential terms
correspond to the continuous diffusion component of the asset price dynamics. 

Now, we evaluate its numerical solution. The resulting PIDE is challenging to solve directly due to the integral term that results from the jump component. In order to solve this, we use the FFT, which converts the equation into Fourier space and effectively handles the nonlocal integral component. We changed the coordinate system in which we solve the problem rather than the financial model itself to make it more robust and computationally feasible. A new physical coordinate $x = log(S)$ and a new solution function $u(x,t) = V(e^x, t)$ are defined, along with a logarithmic transformation applied to the price of the asset. Importantly, it turns multiplicative jumps into additive shifts from $x \rightarrow x+y$, which turns the jump integral into a conventional convolution, and further simplifies the differential operator.
Then the jump integral is written as
\[
\lambda \int_{-\infty}^{\infty} u(x + y, t) f_Y(y) \, dy,
\]
as we established, the log-space transformation conveniently reveals this integral to be a convolution of the solution $u(x, t)$ with the log-jump size probability density function $f_Y(y)$. Further, we employ the Convolution Theorem and a fundamental principle of Fourier analysis, which allows us to replace the complex integration with a sequence of highly optimized FFT operations. Denoting the forward Fourier Transform as $\mathcal{F}$ and the inverse as $\mathcal{F}^{-1}$, the convolution integral is given as:
\[
(u * f_Y)(x) = \mathcal{F}^{-1}\left\{\mathcal{F}\{u(x, t)\} \odot \mathcal{F}\{f_Y(y)\}\right\},
\]
where $\odot$ represents the element-wise product. Furthermore, Equation \ref{eq:merton_pide} applies the chain rule on the derivatives to give the final PIDE in log-space, and is written as
\begin{equation}
    \frac{\partial u}{\partial t} + \left(r - \frac{1}{2}\sigma^2\right)\frac{\partial u}{\partial x} + \frac{1}{2}\sigma^2\frac{\partial^2 u}{\partial x^2} - (r+\lambda)u + \lambda(u * f_Y)(x) = 0.
    \label{eq:pide_log_space2}
\end{equation}

For the FFT-based approach we use, this log-space formulation's convolutional jump term is optimal. However, if $r \approx\frac{1}{2}\sigma^2$, the  drift term  of constant-coefficient $(r - \frac{1}{2}\sigma^2$) may cause numerical instabilities.
As a result, our approach takes a more resilient hybrid approach. To establish our wavelet basis and calculate the jump integral as an efficient convolution, we employ the log-space coordinate $x = log(S)$. However, the derivatives are substituted into Equation \ref{eq:pide_log_space2} to enforce the PIDE residual itself in the original S-space. Further, we apply the chain rule to connect our log-space wavelet derivatives $u_x, u_{xx}$ to the S-space terms

\[\ V_S = \frac{\partial V}{\partial S} = \frac{u_x}{S} \ \ \text{and} \
\ V_{SS} = \frac{\partial^2 V}{\partial S^2} = \frac{u_{xx} - u_x}{S^2}.\]
This approach provides the numerical stability of the S-space PIDE while maintaining the FFT efficiency of the log-space transform. We employ HW-PINN to solve this transformed PIDE in the following section.
\section{Hybrid Wavelet-based Physics-Informed Neural Networks}\label{section:4}

\begin{figure}[htbp] 
    \centering
    \includegraphics[width=0.8\linewidth]{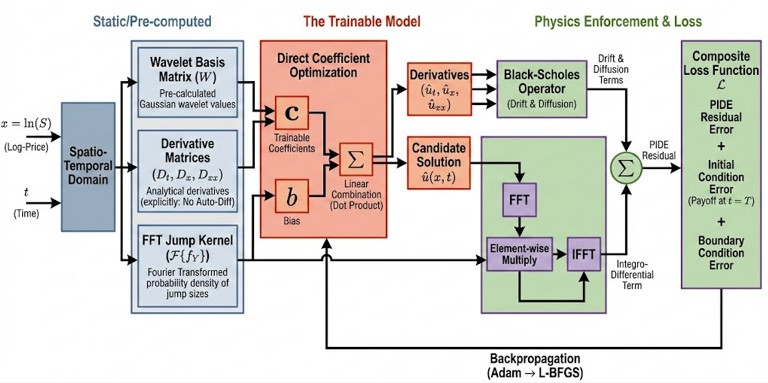} 
    \caption{Architecture diagram of our proposed approach.}
    \label{fig:HW-PINN_architecture}
\end{figure}

Building upon the above problem formulation,  this section describes the systematic framework of HW-PINN to solve the Merton PIDE. Our approach combines a wavelet-based function approximation, an analytical derivative computation, an efficient FFT-based method for calculating the jump integral, and a robust multi-stage optimization approach for training the neural network. Figure \ref{fig:HW-PINN_architecture} depicts the overall architecture of the HW-PINN, which acts as a blueprint for the interrelated parts and how they work together, as we discuss below.

The first stage of our HW-PINN framework is to define the structure of our approximation solution $\hat{u}(x,t)$. In order to define $\hat{u}(x,t)$ as a linear combination, we employ a wavelet family of basis functions, and it is written as
\begin{equation*}
    \hat{u}(x, t) = \sum_{i=1}^{N} c_i \Phi_i(x, t) + \mathfrak{B} \ ,
\end{equation*}
here, this architecture has a unique position, scale, and shape for each of the wavelet basis functions $\Phi_i(x,t)$. $c_i$ represent wavelet coefficients, \(\mathfrak{B}\) denotes a trainable bias parameter, and $N$ is the total number of resolutions.

We now examine how the characteristics of the selected basis functions underpin the effectiveness of this strategy. A single mother wavelet function $\Phi(x)$ is scaled and shifted to create the family of basis functions. In this study, we employ the product of  two one-dimensional Gaussian functions, one for the temporal dimension $(t)$ and the other for the spatial dimension $(x)$, to create a two-dimensional wavelet, which is written as:
\begin{equation*}
\Phi_{j_x, k_x, j_t, k_t}(x, t) = \Phi(X_{\text{scaled}})\Phi(T_{\text{scaled}}) = (X_{\text{scaled}} T_{\text{scaled}}) \exp\left(-\frac{X_{\text{scaled}}^2 + T_{\text{scaled}}^2}{2}\right),
\end{equation*}
where the parameters $j_x, j_t$ control the scale and $k_x, k_t$ control the position, whereas $X_{\text{scaled}} = j_x x - k_x \quad \text{and} \quad T_{\text{scaled}} = j_t t - k_t.$ We ensure that our foundation includes functions that capture both broad, low-frequency market trends and the sharp, localized effects of jumps by selecting a wide range of these parameters. This structure allows for a multi-resolution study of the solution space.

Further, the ability to compute its derivatives analytically is an important aspect of selecting an analytical basis function, such as the Gaussian. This enables us to pre-calculate the derivative values of each basis function as well as the value of the basis function itself at each location in our domain. In order to save computation, we store these values in the following matrices: \(\mathbf{W}\), \(\mathbf{D_tW}\), \(\mathbf{D_xW}\), and \(\mathbf{D_{xx}W}\). Intuitive, highly optimized matrix-vector products can be used to calculate the solution and its derivatives using the following matrices:
\begin{equation*}
    \hat{u} = \mathbf{W}\mathbf{c} + b, \quad \hat{u}{_{t}} = \mathbf{D_tW}\mathbf{c}, \quad \hat{u}{_{x}} = \mathbf{D_xW}\mathbf{c}, \quad \hat{u}{_{xx}}= \mathbf{D_{xx}W}\mathbf{c} \ .
\end{equation*}
The effectiveness of the HW-PINN is based on this. A major obstacle in conventional PINN architectures, automatic differentiation (AD) and the computational network complexity that goes along with it are eliminated.

We now have all the elements required to build the complete PIDE residual and specify the loss function that will guide the training of our model by combining the analytical derivatives from the wavelet matrices with highly efficient FFT-based computation of the jump integral, as addressed in the following subsection.
 \subsection{ Training Architecture and Loss Function}

Here, we combine all of the computational and mathematical components to enable learning for the model.   The process begins with the formulation of a composite loss function, which is then trained to reduce our loss function through a carefully organized, multi-phase training approach that guarantees reliable and effective convergence to an accurate solution.
\subsubsection{Training Data Generation}
In order to enforce the governing equations and boundary constraints via the composite loss function, our approach requires sampling different point sets inside the computational domain.

\begin{enumerate}
    \item \textbf{Collocation Points ($\mathcal{Y}_c$)}: To enforce the PIDE residual $(\mathcal{L}_{PIDE})$, $N_c$ = 8192 collocation points $(x_c, t_c)$ were sampled within the interior of the log-price domain $[x_{min}, x_{max}] \times [0,T]$. In order to provide a more uniform distribution than conventional pseudo-random techniques, a Sobol quasi-random sequence was used for this sampling, which is advantageous for effectively covering the solution space.
    
     \item \textbf{Terminal Condition Points $(\mathcal{Y}_i)$}: In order to strongly enforce the known payoff of the option at expiry  $(\mathcal{L}_{IC})$, $N_ic$ = 1024 terminal condition points ($x_i,T)$ are sampled uniformly at random along the log-price axis exclusively at the final time $t = T$. At these points, the target value is the payoff function $V_{ic} = max (e^{x_i}-K,0)$.
    \item \textbf{Boundary Condition Points ($\mathcal{Y}_b$}): To constrain the spatial boundaries $(\mathcal{L}_{BC})$, $N_{bc}$ = 512 boundary condition points are sampled uniformly at random along the time axis at both the lower boundary $x_{min} = log(S_{min})$ and the upper boundary $x_{max} = log(s_{max})$. The enforced boundary values were $u(x_{min}, t)=0$ and $u(x_{max}, t ) = e^{x_{max}}-Ke^{-r(T-t)}$.
    \subsubsection{Loss Function Formulation}
We now define the total loss function $\mathcal{L}$ of our framework. It consists of three different constraints: the PIDE residual loss, the boundary condition loss, and the initial condition loss. Every constraint in the total loss function makes sure that the wavelet expansion that is learned fits an aspect of the computing problem.

In order to determine how well our wavelet solution fulfills the governing S-space PIDE (see, Equation \ref{eq:pide_log_space2}), we first begin with the PIDE Residual Loss, $\mathcal{L}_{PIDE}$. The residual, as determined by the chain rule derivatives specified in Section \ref{section:3}, is:
\[
\mathcal{L}_{PIDE} = \frac{1}{N_c} \sum_{i=1}^{N_c} \left| \left( \frac{\partial u}{\partial t} + \frac{1}{2} \sigma^2 S^2 V_{SS} + r S V_S - (r + \lambda) u + \lambda (u * f_Y)(x) \right) \right|^2.
\]
The second is the initial condition loss $\mathcal{L}_{IC}$, which ensures that, for a set of $N_{ic}$ points, the solution computes the option's known, non-smooth payoff structure at maturity (t=T).
\begin{equation*}
\mathcal{L}_{IC} = \frac{1}{N_{ic}} \sum_{i=1}^{N_{ic}} \left| \hat{u}(x_i, T) - \max(e^{x_i} - K, 0) \right|^2.
\end{equation*}
The third is the Boundary Condition Loss, $\mathcal{L}_{BC}$, which enforces the behavior of the solution at $x_{upper}$ and $x_{lower}$, the spatial boundaries of the computational domain, over $N_{bc}$ points in time.
\begin{equation*}
\mathcal{L}_{BC} = \frac{1}{N_{bc}} \sum_{i=1}^{N_{bc}} \left( \left| \hat{u}(x_{\text{lower}}, t_i) - V_{\text{left}} \right|^2 + \left| \hat{u}(x_{\text{upper}}, t_i) - V_{\text{right}} \right|^2 \right).
\end{equation*}
The total loss of HW-PINN is the weighted sum of the above losses, which is
\begin{equation*}
\mathcal{L} = \lambda_{PIDE} \mathcal{L}_{PIDE} +  \lambda_{IC} \mathcal{L}_{IC} +  \lambda_{BC} \mathcal{L}_{BC}  .
\end{equation*}
The loss weights $\lambda_{PIDE}, \ \lambda_{IC},\ \lambda_{BC}$ 
 directly determine the impact of each loss on the parameter optimization, which consequently influences our HW-PINN convergence.

\subsubsection{Multi-Stage Training Strategy}
Here, we use a three-phase training technique that capitalizes on the advantages of various optimization algorithms for reliable and effective convergence in order to successfully overcome the intricate and high-dimensional loss environment.

\begin{enumerate}
\item \textbf{Initial Adam Optimization}: In the first stage, our HW-PINN model is trained over a significant number of epochs using the Adam optimizer. Adam is a gradient-based, first-order algorithm that can change its learning rate, making it ideal for quickly searching the parameter space for a potentially lucrative basin of attraction. If the loss stagnates, a ReduceLROnPlateau learning rate scheduler dynamically lowers the learning rate to improve convergence.

\item \textbf{Adam Refinement Stage}: In the second stage, a Coefficient Refinement Network is initialized using the bias and best-performing wavelet coefficients from the first training. The Adam optimizer is then used to train the network for additional epochs (Refinement Epochs) at a much lower learning rate (Refinement\_ADAM\_LR). Additionally, this step uses gradient clipping to ensure more stable, fine-grained convergence within the target basin by preventing overly large parameter updates.

  \item \textbf{L-BFGS Finalization}: The final stage uses the L-BFGS optimizer. These quasi-Newton techniques use second-order information to achieve faster, more precise convergence to a sharp local minimum. It is invoked after the Adam refinement stages to perform a final, high-precision tuning of the wavelet coefficients. This practice of using L-BFGS as the final optimizer is common and effective in our HW-PINN.
    
\end{enumerate}
\end{enumerate}

We have discussed the computation of our HW-PINN in Algorithm \ref{alg:HW-PINN}.

\renewcommand{\algorithmicrequire}{\textbf{Input:}}
\renewcommand{\algorithmicensure}{\textbf{Output:}}
\begin{algorithm}[H]
\caption{Computation of HW-PINN for Merton PIDE}
\label{alg:HW-PINN}
\begin{algorithmic}[1]
\REQUIRE

    \parbox[t]{\dimexpr\linewidth-\algorithmicindent-\algorithmicindent}{ 
    Model parameters $\mathcal{P} = \{r, \sigma, K, T, \lambda, \mu_J, \sigma_J\}$; \\
    Domain $S \in [S_{\min}, S_{\max}]$, $t \in [0,T]$; \\
    Wavelet scale ranges $(j_x^{\min}, j_x^{\max})$, $(j_t^{\min}, j_t^{\max})$; \\
    Training sizes $N_c, N_i, N_b$; weights $w_p, w_i, w_b$
    }

\ Optimized coefficients $\mathbf{c}^*$, bias $b^*$

\ENSURE Transform to log-price: $x = \log S$

\STATE \parbox[t]{\dimexpr\linewidth-\algorithmicindent}{ 
    Generate Gaussian wavelet basis $\Phi_i(x,t)$ for $i=1...N_f$ using parameters
    $(j_x, j_t, k_x, k_t)$ from wavelet family $\mathcal{F}$:
    \[
    \begin{split} 
    \Phi_i(x,t) = (j_x x - k_x)(j_t t - k_t) \\
     \times e^{-\frac{1}{2}\left[(j_x x - k_x)^2 + (j_t t - k_t)^2\right]}
    \end{split}
    \]
    over scales $j_x, j_t$ and shifts $k_x, k_t$ within domain.
    } 

\STATE Sample collocation points $(x_c, t_c) \sim \mathrm{Sobol}(N_c)$, initial $(x_i, T)$, boundary $(x_{\min}, t_b), (x_{\max}, t_b)$

\STATE Precompute FFT kernel $\hat{f}_Y = \mathcal{F}[\mathrm{IFFTShift}(f_Y)]$ with $f_Y \sim \mathcal{N}(\mu_J, \sigma_J)$

\STATE Initialize $\mathbf{c} \sim \mathrm{XavierUniform}(\mathbb{R}^N)$, $b = 0.5$

\WHILE{not converged}
    
\STATE \parbox[t]{\dimexpr\linewidth-\algorithmicindent}{ 
        Compute loss $\mathcal{L}_{\text{total}} =  \lambda_{PIDE} \mathcal{L}_{PIDE} +  \lambda_{IC} \mathcal{L}_{IC} +  \lambda_{BC} \mathcal{L}_{BC}$ where:
        \begin{align*}
            \mathcal{L}_{\text{PDE}} &= \mathbb{E}_{(x,t)}\left[\mathcal{R}_{\text{PIDE}}(x,t)^2\right], \\
            \mathcal{L}_{\text{IC}} &= \mathbb{E}_{x_i}\left[(u(x_i,T) - (e^{x_i} - K)^+)^2\right], \\
            \mathcal{L}_{\text{BC}} &= \mathbb{E}_{t_b}\left[u(x_{\min},t_b)^2 + (u(x_{\max},t_b) - (e^{x_{\max}} - Ke^{-r(T-t_b)}))^2\right],
        \end{align*}
        and the residual in physical space is:
        \[
        \mathcal{R}_{\text{PIDE}} = u_t + \frac{1}{2}\sigma^2 S^2 V_{SS} + r S V_S - (r+\lambda)u + \lambda \int u(x+y,t) f_Y(y)\,dy,
        \]
        with $S = e^x$, $V_S = u_x / S$, $V_{SS} = (u_{xx} - u_x)/S^2$, computed via analytical derivatives of the wavelet expansion:
        \[
        \hat{u}(x, t) = \sum_{i=1}^N c_i \Phi_i(x,t) + \mathfrak{B}.
        \]
        } 

\STATE Evaluate integral term using FFT-based convolution.
    
\STATE \parbox[t]{\dimexpr\linewidth-\algorithmicindent}{ 
        Update $(\mathbf{c}, b)$ using multi-stage optimization:
        \begin{itemize} 
            \setlength\itemsep{0em} 
            \item Stage I: Adam ($\text{lr} = \text{LR}_1$)
            \item Stage II: Adam with gradient clipping ($\text{lr} = \text{LR}_2 < \text{LR}_1$)
            \item Stage III: L-BFGS (with strong-Wolfe line search)
        \end{itemize}
        } 
\ENDWHILE

\RETURN $\mathbf{c}^*, b^*$
\end{algorithmic}
\end{algorithm}

\section{Experimental Result and Discussion}\label{section:5}
In this section, we present the empirical results of our study, validating the accuracy and robustness of our proposed HW-PINN model. We demonstrate its performance against an established benchmark across six realistic market scenarios, providing a clear picture of the model's capabilities as a practical financial pricing tool. We will pay special attention to the Low Jump Intensity case to showcase the model's strong performance and provide a detailed analysis of what these results mean in both technical and financial contexts.

\begin{figure}[htbp] 
    \centering
    \includegraphics[width=0.8\linewidth]{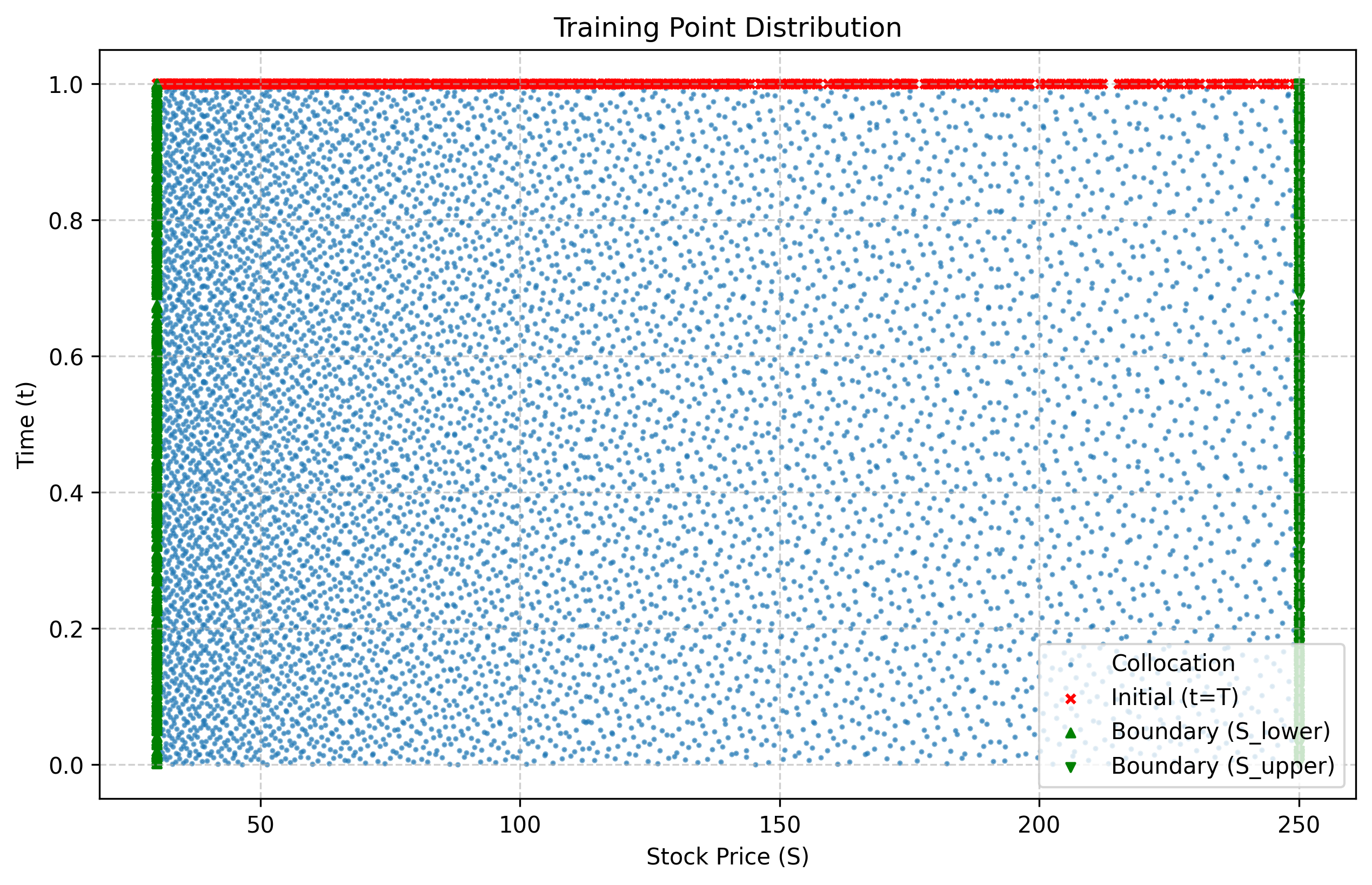} 
    \caption{Distribution of the 8,192 training points (blue), initial condition points (red), and boundary condition points (green) across the price-time domain.}
    \label{Fig:2}
\end{figure}

 All numerical experiments are conducted on an HPC cluster equipped with NVIDIA H100 NVL GPUs, which significantly reduces training time and accelerates computation. This configuration provides consistent GPU acceleration while providing efficient resource sharing across large computational workloads.
 The implementation is carried out in Python using the PyTorch deep learning framework. 
The loss function weights are held constant at $\lambda_{PIDE}=1.0$, $\lambda_{IC}=5000.0$, and $\lambda_{BC}=10.0$ to strongly enforce the initial and boundary conditions. To assess the model's performance, we compare its predictions against a high-fidelity Carr-Madan FFT-based solver  \cite{carr1999option}. We tested the model across six distinct market scenarios, each designed to probe its robustness under different financial conditions. The HW-PINN model is constructed with a wavelet family of over 2,500 basis functions, and the loss function is evaluated on a set of 8,192 collocation points, quasi-randomly distributed across the domain as depicted in Figure \ref{Fig:2}.

\subsection{Model Performance across Several Market Scenarios}
The primary finding of our study is that the HW-PINN model produces extremely accurate option prices across a variety of normal and challenging market conditions. The relative percentage error for four at-the-money and in-the-money validation points across the six market scenarios is displayed in Table \ref{tab:1}, which summarizes the findings.  Table \ref{tab:1} shows that the error is typically less than 2\% in a majority of cases, which is a great outcome for a numerical PDE solver. The model's potential utility for real-world pricing applications is demonstrated by its success under challenging but realistic market regimes.
Our model performs exceptionally well in smooth market conditions, achieving very high accuracy in the low jump intensity case, with errors constantly below 1\%. The reliability of our HW-PINN in managing standard and highly volatile markets is demonstrated by the errors remaining reasonably low (usually within 1–5\%) under the baseline and high volatility scenarios. The model maintains an acceptable level of accuracy in challenging situations, such as high interest rates, crash scenarios, and high jump intensity; nevertheless, in some situations, significantly higher inaccuracies (up to about 10\%) are noted because of increasing nonlinearity and discontinuities. Overall, the outcomes demonstrate our HW-PINN's ability to simulate both smooth and irregular behaviors governed by PIDEs and show that it effectively captures challenging option pricing dynamics across an extensive range of financial situations.
However, in a variety of market scenarios, Table \ref{tab:1} demonstrates that PINN produces large errors while failing severely. This obstacle comes from the Merton PIDE's jump integral term, which introduces non-local behavior that PINNs find difficult to capture.
 However, as shown in Table \ref{tab:1}, HW-PINN maintains stability and produces superior results than PINN and Benchmark because its wavelet basis structure regularizes the solution space, making it more resistant for complex jump-diffusion conditions.

\subsubsection{Low Jump Intensity Case}
\begin{figure}[htbp] 
    \centering
    \includegraphics[width=1.2\linewidth]{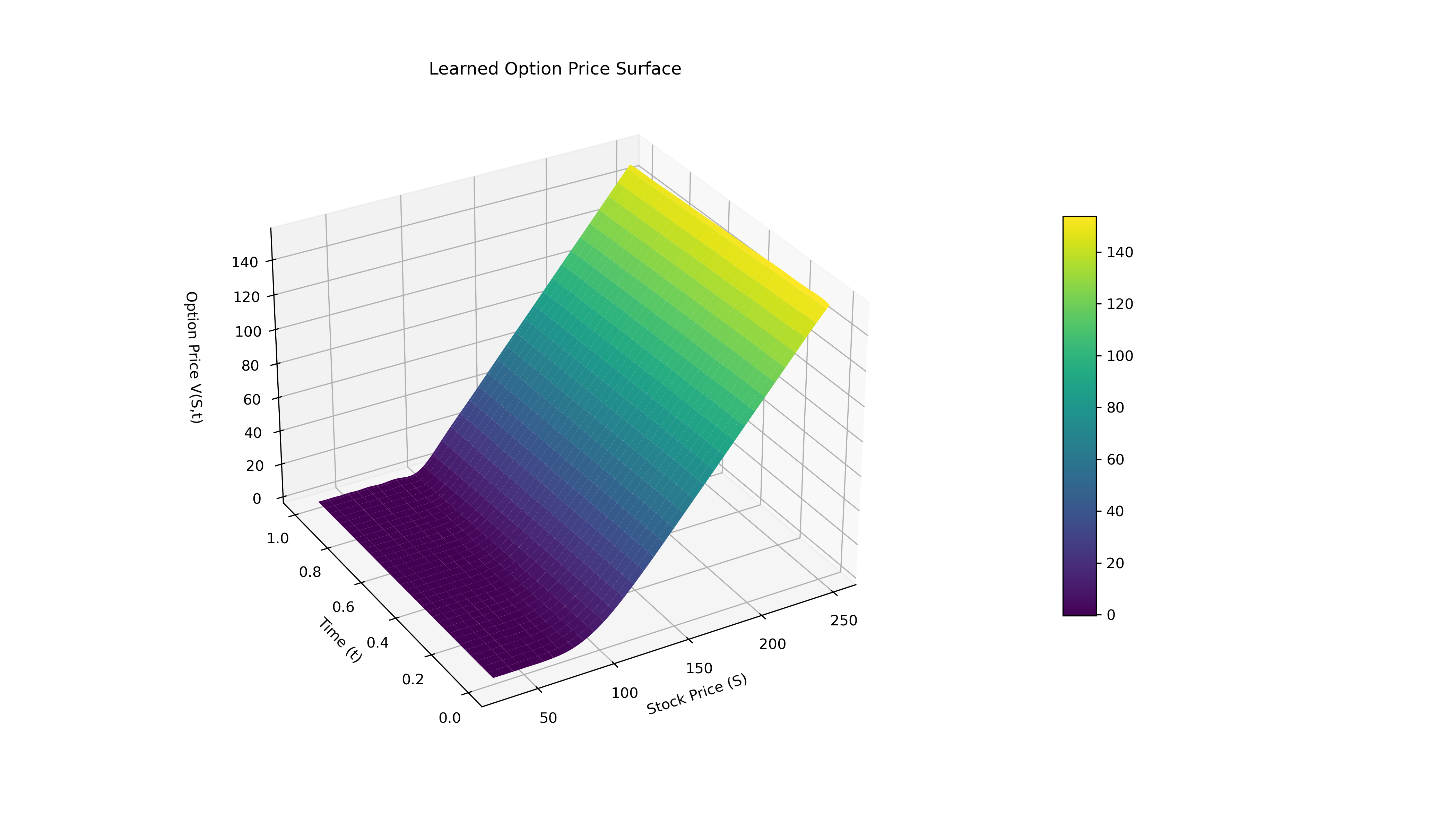} 
    \caption{ HW-PINN learns the option price surface.}
    \label{fig:2}
\end{figure}

\begin{figure}[htbp] 
    \centering
    \includegraphics[width=0.8\linewidth]{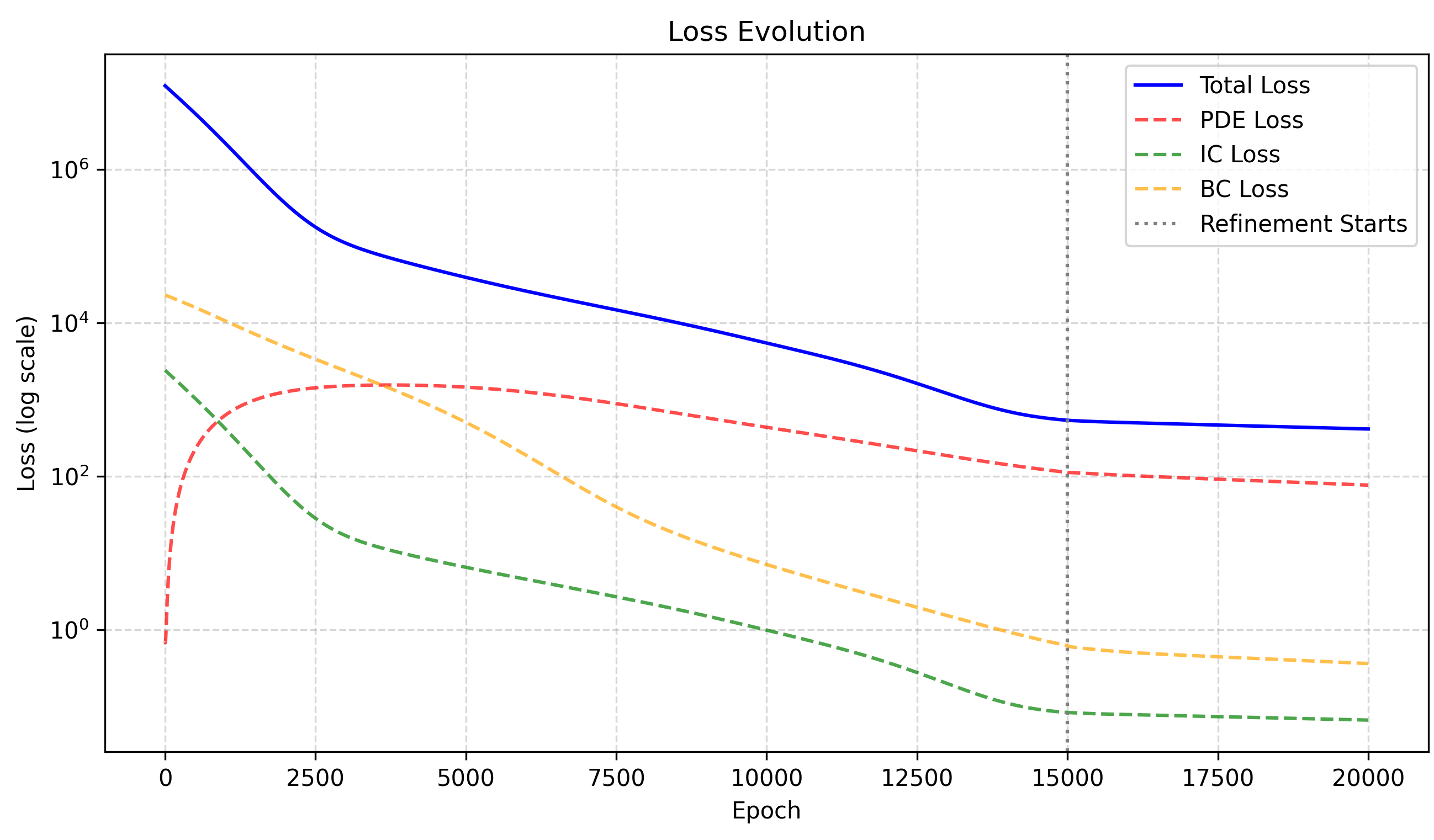} 
    \caption{Training characteristics in the low jump intensity scenario.}
    \label{fig:4}
\end{figure}

\begin{figure}[htbp]
    \centering
    \includegraphics[width=0.48\textwidth]{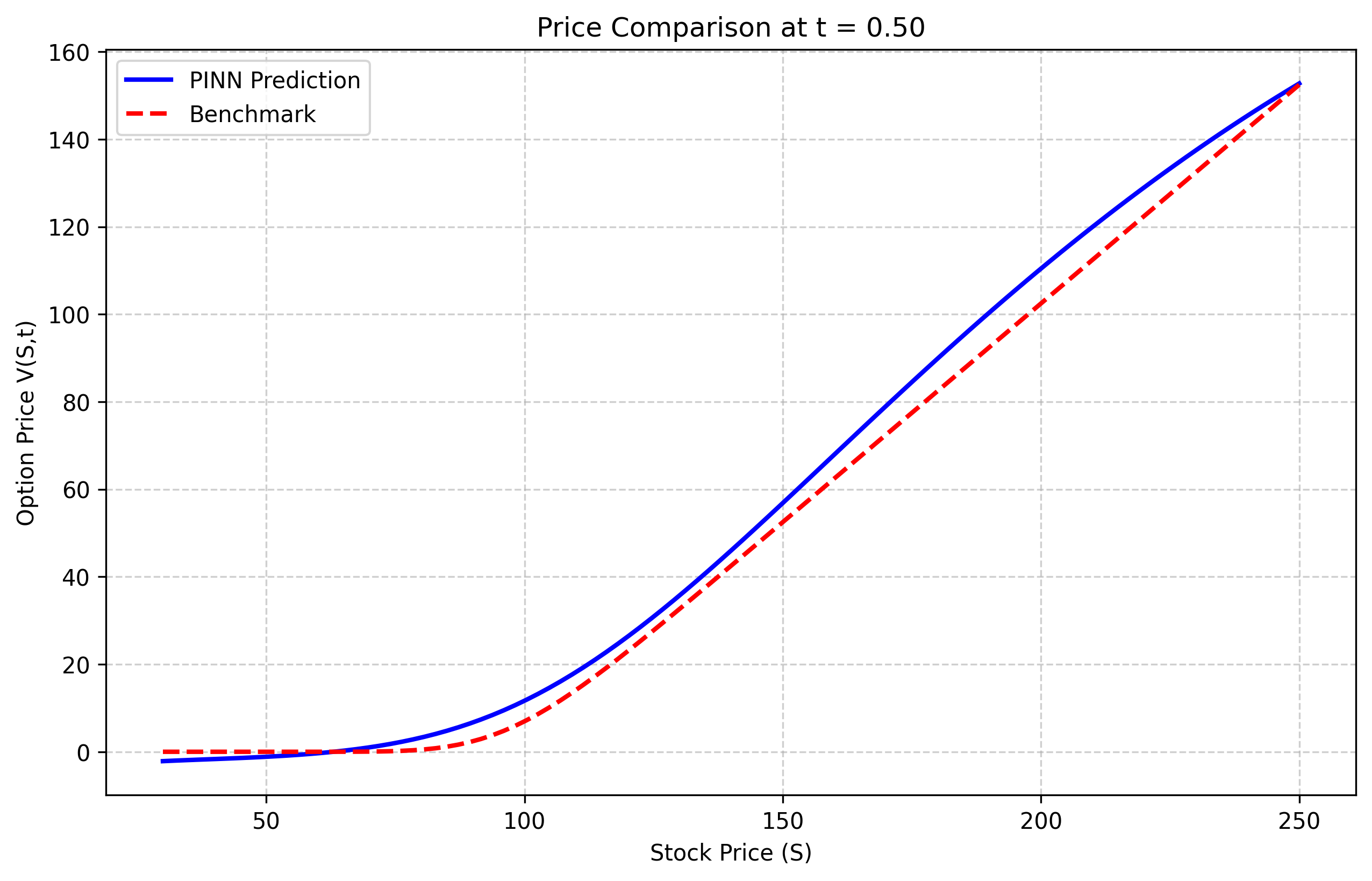}
    \hfill
    \includegraphics[width=0.48\textwidth]{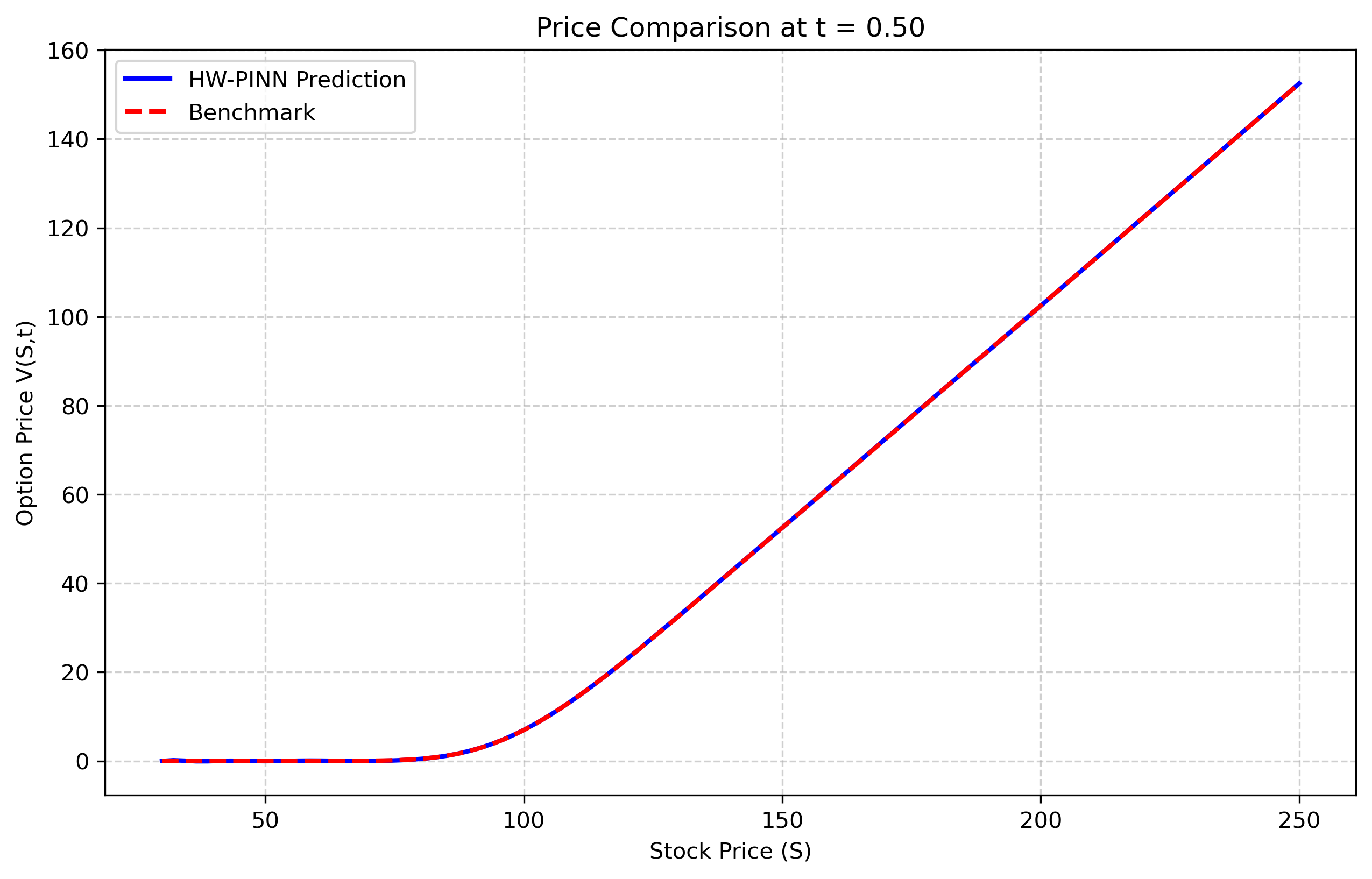}
    \caption{Comparison of PINN and HW-PINN with benchmark solution.}
      \label{fig:3}
\end{figure}

The performance of our model is particularly exceptional in the low-jump-intensity case, which reflects a typical, well-behaved market. In this case, the model achieved a mean relative error of only 0.109\%. To analyze our success, we can examine the entire option-pricing surface generated by the HW-PINN, as shown in Figure \ref{fig:2}. The surface is smooth, monotonic, and convex, adhering to all the theoretical requirements of a no-arbitrage option pricing model. Additionally, Figure \ref{fig:4}, which shows the model's training diagnostics, reveals the underlying technology behind this result.  A stable, effectively convergent loss function with a monotonic decrease is depicted in Figure \ref{fig:4}. This is the most distinctive aspect of a robust and appropriately converged optimization process, showing that the accuracy of our model is an outcome of a reliable learning process rather than coincidence.
Furthermore, we take a cross-section of this surface at a specific point in time and compare it to the benchmark to evaluate accuracy more directly. 
Compared to HW-PINN and PINN solutions, Figure  \ref{fig:3} shows that the PINN approach fails to accurately depict the dependence of option pricing on the underlying asset, with a relatively deviant pattern across stock prices. This behavior is clarified by learning nonlocal jump dynamics being challenging and an imbalance in loss elements.
 On the other hand, HW-PINN exhibits higher resilience and reliability, continuously producing precise outcomes throughout the entire domain. We can see that the almost perfect match between the HW-PINN's prediction and the benchmark price at t=0.50 is displayed in Figure \ref{fig:3}, which clearly validates the high fidelity of the solution.

\begin{figure}[htbp] 
    \centering
    \includegraphics[width=1.2\linewidth]{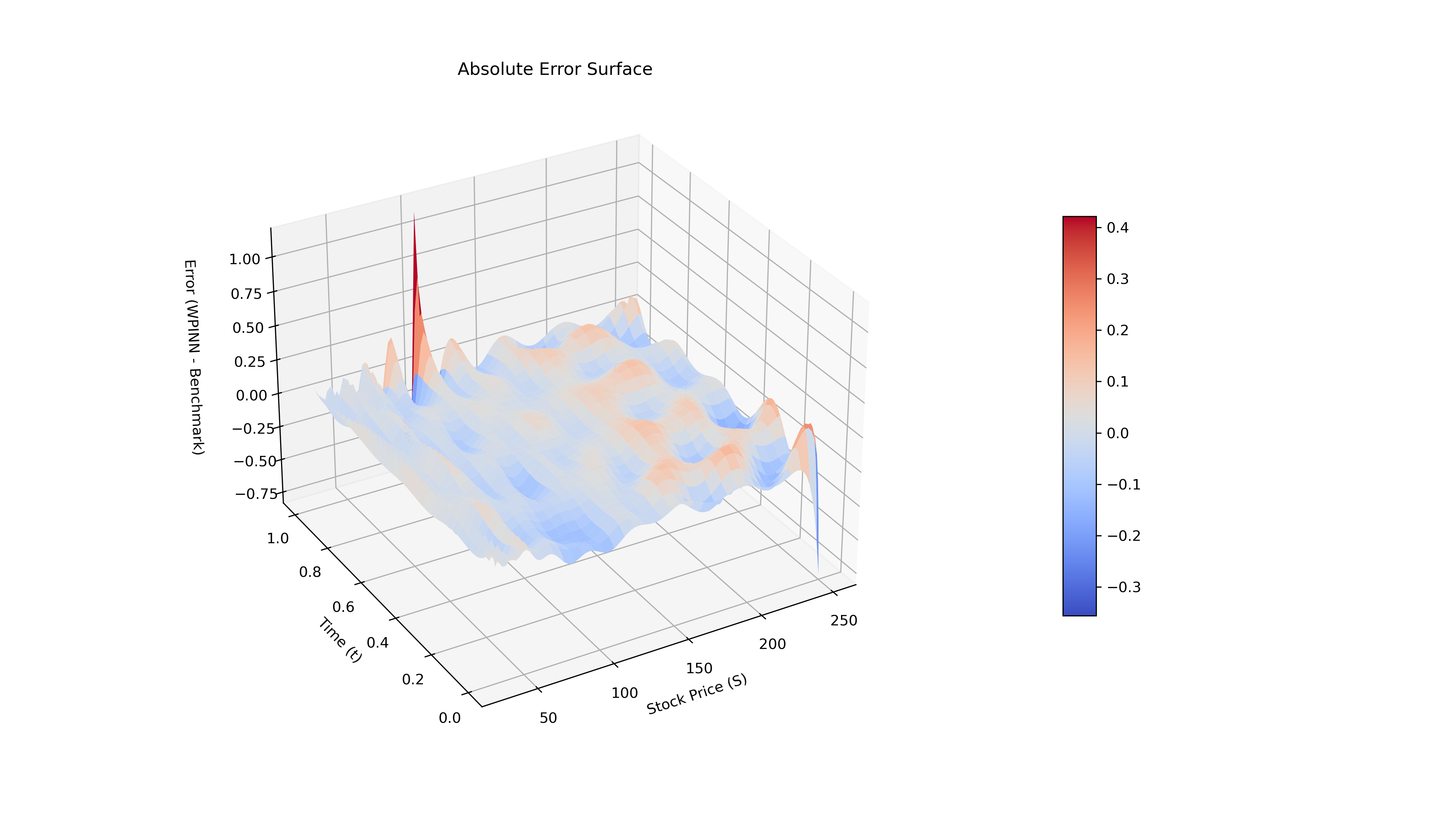} 
    \caption{Comparison prediction of our HW-PINN and the Benchmark solution.}
    \label{4_3D_Absolute_Error_Surface}
\end{figure}

Furthermore, the 3-D absolute error surface between the benchmark and HW-PINN solution is displayed in Figure \ref{4_3D_Absolute_Error_Surface}. The error's sign and magnitude are shown on the color map, where blue areas show the HW-PINN model's underestimating and red areas show a small overestimation.  Overall, the surface shows strong agreement with the benchmark solution, with error ranging from $-0.76$ to $0.25$. Due to the strong curvature of the option price function, small deviations occur close to the strike price area ($S \approx 110$). Still, errors are negligible in the deep out-of-the-money and deep in-the-money regions, where the option value behaves more smoothly. The HW-PINN framework gives a reliable and accurate approximation of the option pricing solution, as seen by the lack of significant oscillations across time and stock price.

\begin{figure}[htbp] 
    \centering
    \includegraphics[width=0.8\linewidth]{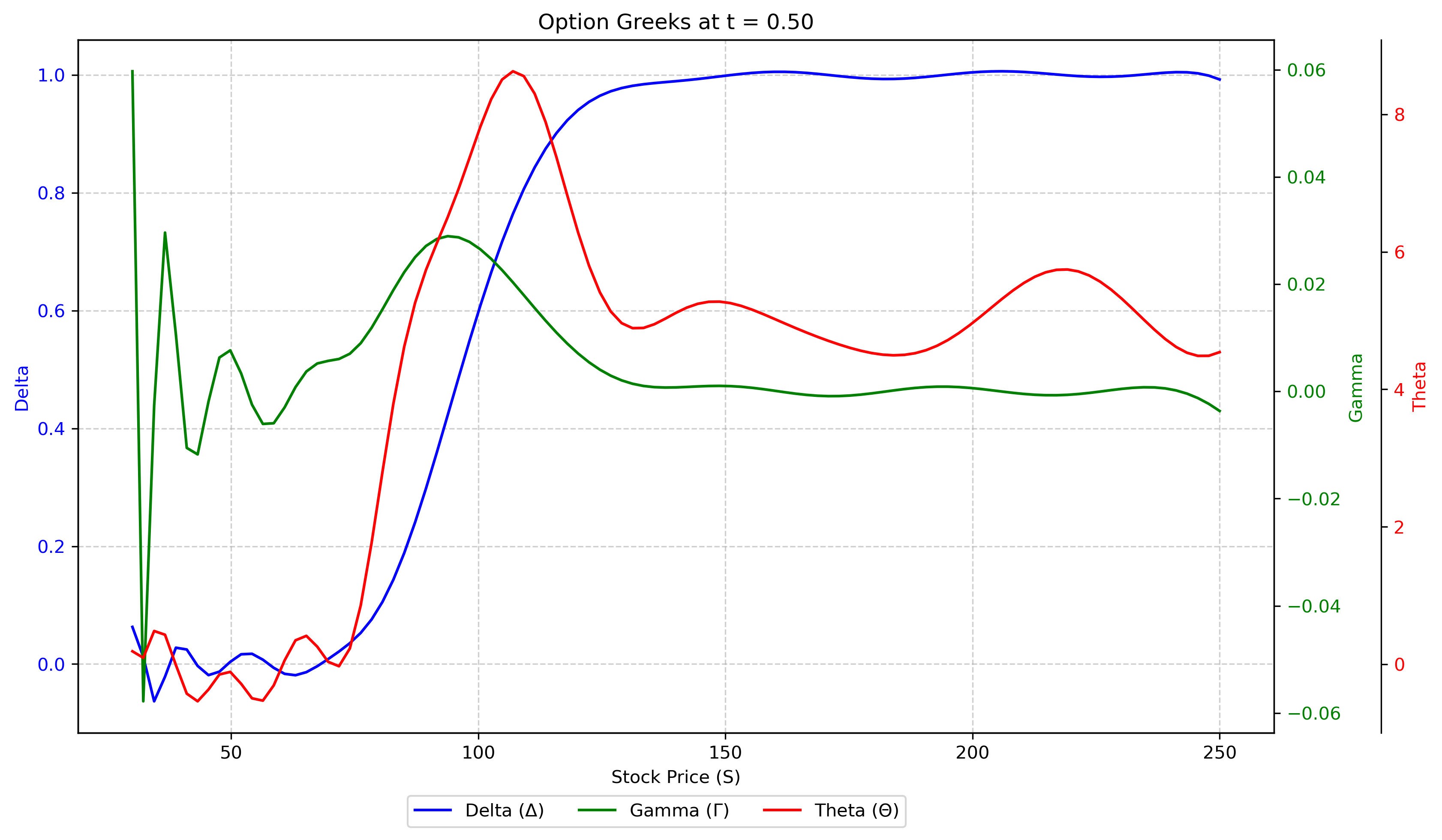} 
    \caption{Option Greeks for the low jump intensity case.}
    \label{fig:5}
\end{figure}

Finally, the option Greeks are highly concentrated in the at-the-money (ATM) area under a low jump intensity, as depicted in Figure \ref{fig:5}. Figure \ref{fig:5} shows that Delta ($\Delta$) is smooth and stays near zero in the deep out-of-the-money (OTM) region, indicating the option value's minimal sensitivity to fluctuations in the price of the underlying asset. Delta rises quickly as the stock price approaches the ATM, indicating a dramatic shift in price sensitivity. Further, Delta progressively approaches unity in the deep in-the-money (ITM) region, indicating that the option price moves nearly in parallel with the underlying asset. Between the OTM and ITM regions, Gamma ($\Gamma$) decreases toward zero and reaches its maximum near the ATM region. This sharp peak indicates that the option price curve with respect to the underlying asset is concentrated around the target price, where Delta changes most rapidly. The option price becomes roughly linear in these scenarios, as evidenced by the disappearance of Gamma in the most extreme regions. Although  Theta ( $\Theta$) achieves its maximum magnitude close to the ATM region, indicating that time decay is most important when the option has the highest time value. Distinct from the strike price, Theta significantly decreases, suggesting that for options that are ITM or OTM, the effect of time decay becomes insignificant. The behavior of $\Delta$, $\Gamma$, and $\Theta$ together shows that option sensitivities are heavily concentrated around the strike price, given the market conditions under consideration. This localization of sensitivities suggests that positions far from the strike show comparatively steady and reliable behavior, whereas risk exposure and hedging requirements are most pronounced in the ATM region.

\subsection{Risk Analysis Using Value at Risk and Conditional Value at Risk}
In this subsection, we evaluate the option price accuracy of our proposed HW-PINN approach and analyze the financial risk associated with the option prices under different market scenarios.

\begin{table}[ht!]
\centering
 \renewcommand{\arraystretch}{5}
\rowcolors{2}{blue!10}{blue!10} 
    \renewcommand{\arraystretch}{1.2} 
\begin{tabular}{|p{4cm}|p{.7cm}|p{.8cm}|p{1.5cm}|p{1.7cm}|p{2.2cm}|p{1.6cm}|p{1.5cm}|}
\rowcolor{green!30} 
\toprule
\textbf{Test Case} & \textbf{S} & \textbf{t} & \textbf{PINN} & \textbf{HW-PINN} & \textbf{Benchmark \cite{carr1999option}} & \textbf{HW-PINN Error (\%)} & \textbf{PINN Error (\%)} \\
\midrule
Low Jump Intensity & 90 & 0.75 & 2.9206 & 0.9165 & 0.9227 & \textbf{0.672} & 216.530  \\
                   & 100 & 0.50 & 6.6714  & 6.9643 & 6.9558 & \textbf{0.123} & 4.088 \\
                   & 110 & 0.25 &  13.51852 & 15.9959 & 16.0271 & \textbf{0.195} & 15.652 \\
                   & 150 & 0.75 & 52.6276 & 51.1913 & 51.2440 & \textbf{0.103} & 2.700  \\
\midrule
Baseline          & 90 & 0.75 &  2.6304 & 1.1076 & 1.0485 & \textbf{5.637} & 150.886  \\
                  & 100 & 0.50 & 6.4053 & 7.2097 & 7.2845 & \textbf{1.027} & 12.069 \\
                  & 110 & 0.25 &  13.8423  & 16.1337 & 16.4175 & \textbf{15.686} &1123.3 \\
                  & 150 & 0.75 &   51.8157  & 51.3479 & 51.2538 & \textbf{1.096} & 427.8 \\
\midrule
High Volatility   & 90 & 0.75 & 3.6714         & 4.3266 & 4.3458 & \textbf{0.441} & 15.517\\
                  & 100 & 0.50 &  8.9238       & 12.8120 & 13.1217 & \textbf{2.361} & 31.992   \\
                  & 110 & 0.25 & 18.2477          & 22.3667 & 22.9828 & \textbf{2.681} & 20.603\\
                  & 150 & 0.75 &  51.6623        & 51.6217 & 51.5957 & \textbf{0.050} & 0.129   \\
\midrule
High Interest Rate & 90 & 0.75 &  5.3713 & 1.4105 & 1.2814 & \textbf{10.078} & 2631.4 \\
                   & 100 & 0.50 & 10.8848  & 8.6683 & 8.6676 & \textbf{0.008} & 303.8 \\
                   & 110 & 0.25 & 19.9935 & 18.9917 & 19.1414 & \textbf{0.782}& 82.8 \\
                   & 150 & 0.75 & 53.8725 & 52.6781 & 52.4787 & \textbf{0.380}& 33.3 \\
\midrule
Crash Scenario    & 90 & 0.75 &  3.0433 & 1.0389 & 1.0910 & \textbf{4.779} & 178.942  \\
                  & 100 & 0.50 &  5.8017 & 7.2321 & 7.7256 & \textbf{6.388} & 24.902\\
                  & 110 & 0.25 & 15.3264  & 16.2838 & 17.1051 & \textbf{4.801} & 10.398 \\
                  & 150 & 0.75 & 54.075  & 51.3880 & 51.3196 & \textbf{0.133} & 5.371  \\
\midrule
High Jump Intensity & 90 & 0.75 &  3.6377 & 1.4969 & 1.4012 & \textbf{6.832} & 159.612 \\
                    & 100 & 0.50 &  6.8979  & 7.7593 & 8.1539 &  \textbf{4.840} & 15.403   \\
                    & 110 & 0.25 &  10.6425  & 16.5121 & 17.4365 &  \textbf{21.301} & 38.964 \\
                    & 150 & 0.75 &  53.6430  & 51.6824 & 51.2902 & \textbf{0.765} & 4.587  \\
\bottomrule
\end{tabular}
\caption{Comparison of PINN, HW-PINN and Benchmark performance under various market scenarios.}
\label{tab:1}
\end{table}

 We employ two popular risk measures: Conditional Value-at-Risk (CVaR) and Value-at-Risk (VaR). The numerical values of these risk measures are presented in Table~\ref{tab:2}.
 Table~\ref{tab:2} demonstrates the option price accuracy of the HW-PINN framework along with the estimated VaR and CVaR in different market conditions. The mean relative error assesses the accuracy of the option price model, while VaR and CVaR values provide insight into the potential downside risk associated with the predicted option prices.  However, VaR and CVaR remain low, indicating that the loss distribution is stable. Our model provides very high accuracy with a very small error of 0.2732\% in the low jump intensity scenario, and the VaR and CVaR values are comparatively low, suggesting a stable and low-risk situation. Both VaR and CVaR grow sharply as market uncertainty increases, as in the high volatility case, indicating higher potential losses even though the prediction error stays small. Our model sustains balanced performance with moderate error and risk levels under a baseline scenario. The mean relative error rise in more extreme situations, such as crash conditions and high jump intensity, emphasizing the challenge of capturing abrupt market shifts, while the difference between VaR and CVaR denotes higher tail risk. The model displays moderate error with regulated risk measures for the high interest rate scenario.

\begin{figure}[htbp] 
    \centering
    \includegraphics[width=0.8\linewidth]{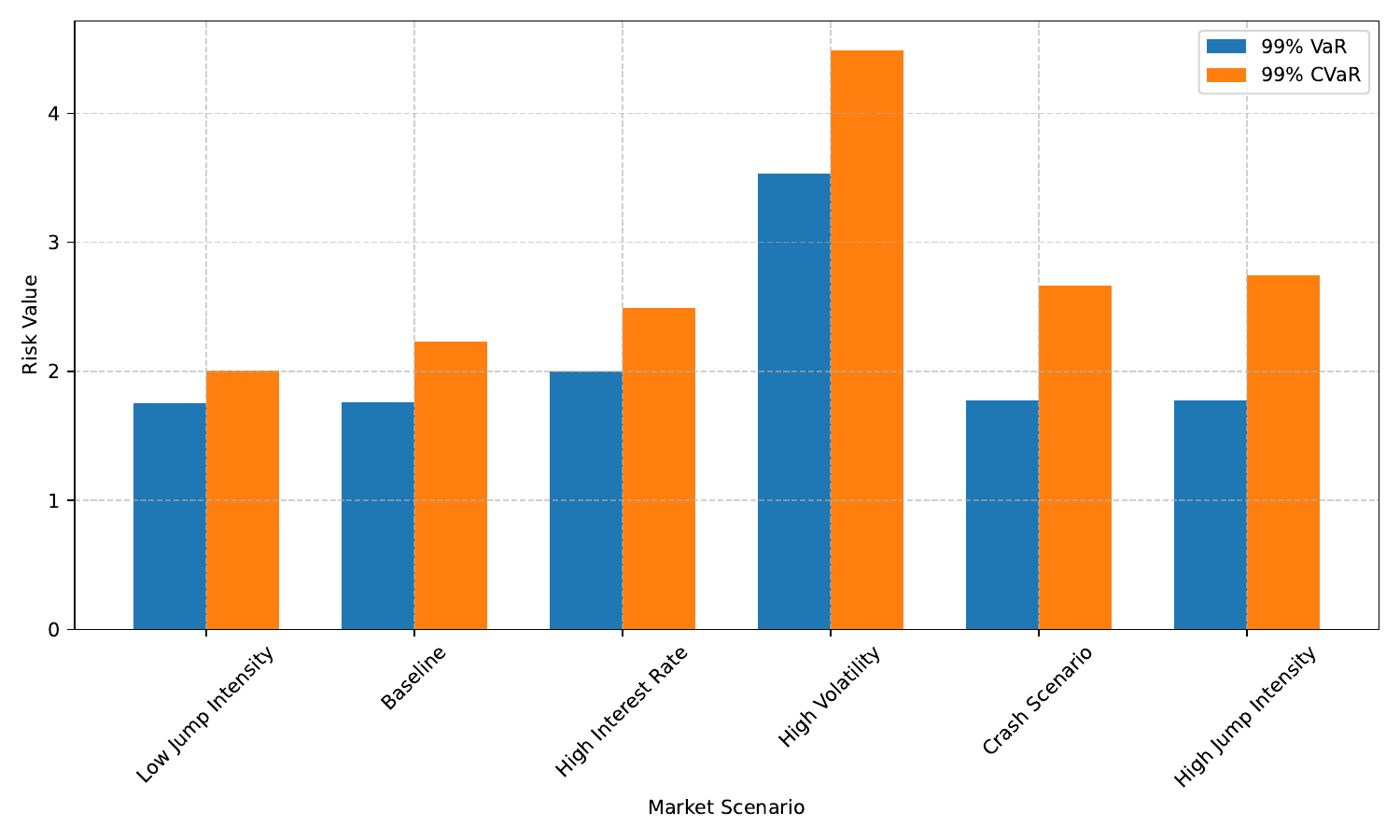} 
    \caption{Comparison of VaR and CVaR across various market scenarios.}
    \label{fig:8}
\end{figure}

\begin{figure}[htbp] 
    \centering
    \includegraphics[width=0.8\linewidth]{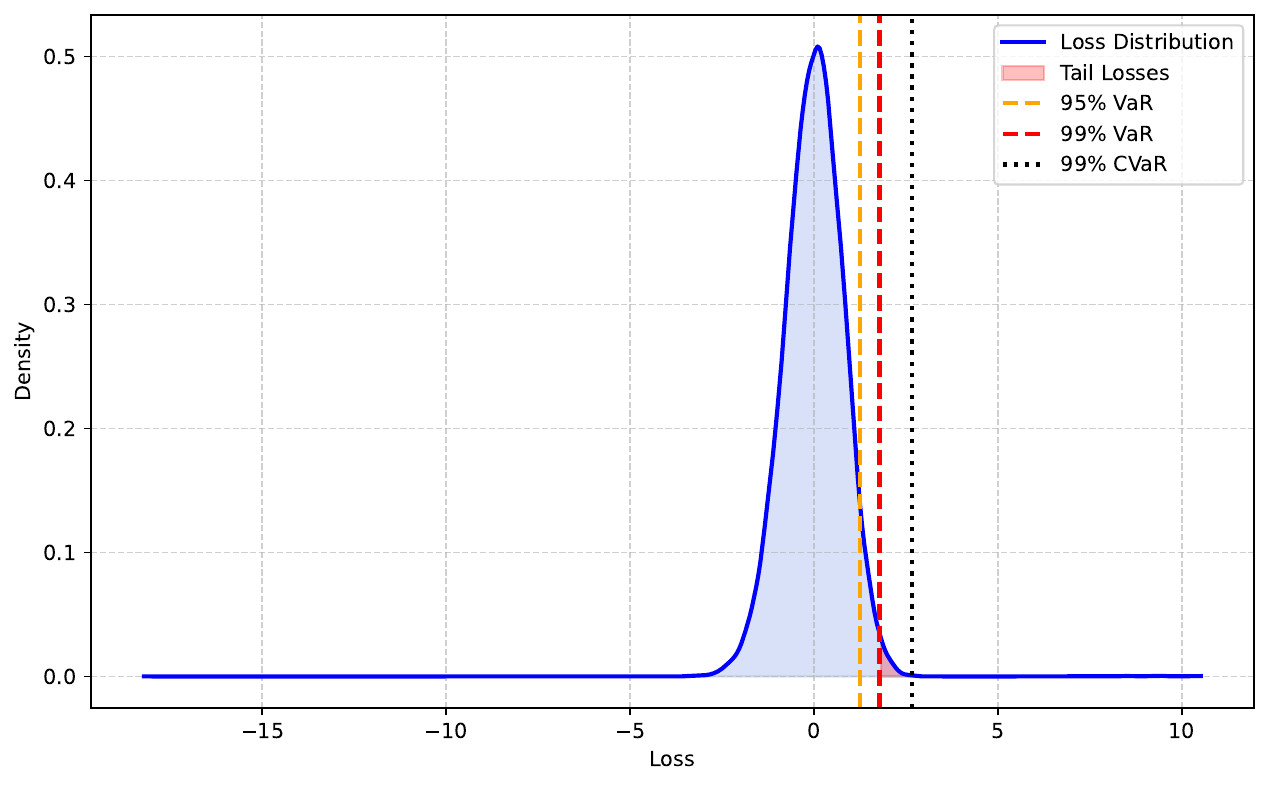} 
    \caption{Loss distribution with VaR and CVaR.}
    \label{fig:9}
\end{figure}

\begin{figure}[htbp] 
    \centering
    \includegraphics[width=0.8\linewidth]{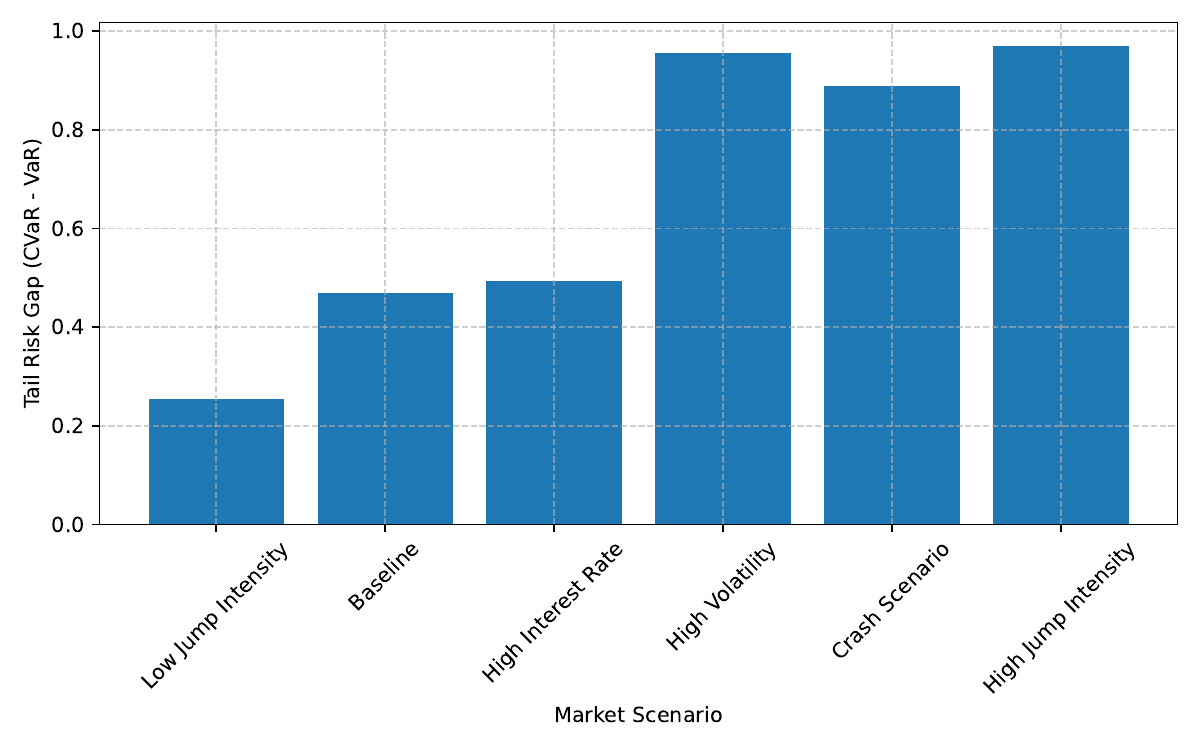} 
    \caption{Tail risk severity across various market scenarios.}
    \label{fig:10}
\end{figure}

To further investigate how the risk measures vary across different market conditions, Figure \ref{fig:8} presents a graphical comparison of the VaR and CVaR values for each scenario.  Figure \ref{fig:8} demonstrates how jump dynamics and volatility affect downside risk and emphasizes the significance of considering both VaR and CVaR into account when evaluating financial risk. It is also important to determine how the underlying loss distribution arises from the underlying loss distribution. Therefore, Figure \ref{fig:9} shows the simulated loss distribution along with the associated VaR and CVaR thresholds. VaR does not accurately reflect the severity of extreme tail scenarios, despite identifying the threshold of possible losses. Figure \ref{fig:10} shows the tail risk gap for each market scenario considered to assess the magnitude of losses exceeding the VaR level. When assessing financial risk in extreme market conditions, it is crucial to consider both VaR and CVaR. Figure \ref{fig:10} demonstrates how volatility and jump intensity have a substantial impact on the severity of tail risk. The combined analysis delivers an in-depth assessment of the proposed HW-PINN approach in several market scenarios. The findings demonstrate that the model not only achieves consistent accuracy in option pricing but also accurately captures the behavior of extreme loss events and downside risk in simulated market scenarios. These results show that our proposed method can deliver significant risk insights and accurate pricing, both of which are critical for financial risk management and option market decision-making.

\begin{table}[ht!]
\centering
\rowcolors{2}{blue!10}{blue!10} 
\begin{tabular}{lccc}
\rowcolor{green!30} 
\hline
\textbf{Case} & \textbf{Mean Relative Error (\%)} & \textbf{99\% VaR} & \textbf{99\% CVaR} \\
\hline
Low Jump Intensity  & 0.2732   & 1.7548 & 2.0085 \\
High Volatility     & 1.3832   & 3.5354 & 4.4903 \\
Baseline            & 2.1440   & 1.7644 & 2.2330 \\
Crash Scenario      & 4.0252   & 1.7775 & 2.6403 \\
High Jump Intensity & 4.4345   & 1.7757 & 2.7444 \\
High Interest Rate  & 2.8120   & 1.9999 & 2.4936 \\
\hline
\end{tabular}
\caption{ Risk analysis of HW-PINN under various market scenarios.}
\label{tab:2}
\end{table}

\section{Conclusion}\label{section:6}
In this study, we present portfolio optimization using the Merton jump-diffusion model for asset prices. 
To accurately reflect market conditions, we integrated European option pricing into traditional portfolio optimization. The Merton jump-diffusion model, which accounts for both continuous market movements and unpredictable jumps in asset prices, is motivated by the observation of abrupt, discontinuous price movements. Ensuring that the option pricing and portfolio optimization challenges within this framework, Merton PIDEs are computationally intractable for conventional numerical methods. We applied the FFT to compute the integral terms arising from the jump component, thereby overcoming these difficulties and reducing computational cost. Further, we used an HW-PINN architecture to enhance the model's ability to capture localized variations and non-smooth behaviors in Merton jump-diffusion models. The wavelet-based multi-resolution structure improves the neural network's capacity to approximate, especially near abrupt gradients and discontinuities. Numerical evaluations show that our proposed approach generates reliable, robust solutions, and the results are in good agreement with established benchmark techniques. We also discussed risk analysis using VaR and CVaR, which highlights the importance of taking both VaR and CVaR into account in risk assessment by revealing the losses under different market situations. This approach may be extended to higher-dimensional portfolio optimization and multi-asset markets in future research.

\section*{Acknowledgements}
Sanjay Kumar Mohanty acknowledges partial support from the Department of Science and Technology, Government of India, under grant SR/FST/MS-II/2023/139(C) at VIT Vellore.

\bibliographystyle{abbrvnat}
\bibliography{references}
\end{document}